\documentclass[reqno]{amsart}

\usepackage{amssymb,epsfig,graphics}
\usepackage{amsmath}
\usepackage{amscd}
\usepackage{verbatim}
\usepackage{youngtab}
\input xypic

\theoremstyle{definition}

\newtheorem{example}{Example}

\def\C{\mathbb C}
\def\R{\mathbb R}
\def\Z{\mathbb Z}

\def\r{\rangle}
\def\l{\langle}

\def\a{\alpha}
\def\o{\omega}

\begin{document}

\title[Polytope ring]
{The rings of $n$-dimensional polytopes}

\author{L. H\'akov\'a}
\address{}
\author{M. Larouche}
\address{}
\author{J. Patera}
\address{Centre de recherches math\'ematiques,
         Universit\'e de Montr\'eal,
         C.P.~6128 Centre-ville,
         Montr\'eal, H3C\,3J7, Qu\'ebec, Canada}
\email{patera@crm.umontreal.ca}

\date{\today}
 \begin{abstract}\

Points of an orbit of a finite Coxeter group $G$, generated by $n$
reflections starting from a single seed point, are considered as vertices of a polytope ($G$-polytope) centered at the origin of a real $n$-dimensional Euclidean space. A
general efficient method is recalled for the geometric description of
$G$-polytopes, their faces of all dimensions and their adjacencies.
Products and symmetrized powers of $G$-polytopes are introduced and
their decomposition into the sums of $G$-polytopes is described.
Several invariants of $G$-polytopes are found, namely the analogs of
Dynkin indices of degrees 2 and 4, anomaly numbers, and congruence
classes of the polytopes. The definitions apply to crystallographic
and non-crystallographic Coxeter groups. Examples and applications
are shown.

 \end{abstract}

\maketitle

\section{Introduction}

Finite groups generated by reflections in a real Euclidean space
$\R^n$ of $n$ dimensions, also called finite Coxeter groups, are split
into two classes: crystallographic and non-crystallographic groups \cite{H,K}.
The crystallographic groups are the Weyl
groups of compact semisimple Lie groups. They are an efficient tool
for uniform description of the semisimple Lie groups/algebras
\cite{Bourbaki,VO,H1}, and they have proven to be an
indispensable tool in extensive computations with the
representations of such Lie groups or Lie algebras (see for example
\cite{GP} and references therein).

Underlying such applications are two facts: (i) Most of the
computation can be performed in integers by working with the weight
systems of the representations involved in a problem, and (ii) the
weight system of a representation of a compact semisimple Lie
group/Lie algebra consists of several Weyl group orbits of the
weights, many of them occurring more than once. Practical importance
of the orbits apparently emerged only in \cite{MP,MP1}, where truly
large scale computations were anticipated.

The crystallographic Coxeter groups are called Weyl groups and
denoted by $W$. Any finite Coxeter group, crystallographic or not,
is denoted by $G$. A difference between the two cases which is of
practical importance to us, is that, lattices with $W$-symmetries
are common crystallographic lattices, while lattices of
non-crystallographic types are dense everywhere in $\R^n$.

Non-crystallographic finite Coxeter groups are of extensive use in
modeling aperiodic point sets with long-range order
(`quasicrystals') \cite{MoPa,RVM,CMP}. Outside traditional
mathematics and mathematical physics, a new line of application of
Coxeter group orbits can be found in \cite{T}; see also the
references therein.

Additional applications of Weyl group orbits are found in
\cite{GPS,ST1,ST2,BS,AP}. Both crystallographic and
non-crystallographic Coxeter groups can be used for building
families of orthogonal polynomials of many variables \cite{D}.

In recent years, another field of applications of $W$-orbits is
emerging in harmonic analysis. Multidimensional Fourier-like
transforms were introduced and are currently being explored in
\cite{P,MoP,KP1,KP2}, where $W$-orbits are used to define families
of special functions, called orbit functions \cite{P}, which serve
as the kernels of the transforms. They differ from the traditional
special functions \cite{KV}. The number of variables, on which the
new functions depend, is equal to the rank of a compact semisimple
Lie group that provides the Weyl group. Two properties of the
transforms stand out: Such special functions are orthogonal when
integrated over a finite region $F$, and they are also orthogonal
when summed up over lattice points $F_M\subset F$. The lattices can
be of any density, their symmetries are prescribed by the Lie
groups. Application of the non-crystallographic groups in Fourier
analysis is at its very beginning \cite{MNP}.

In this paper we have no compelling reason to distinguish
crystallographic and non-crystallographic reflection groups of
finite order. Hence, we consider all finite Coxeter groups
 although from the infinitely many finite Coxeter groups
in 2 dimensions (symmetry groups of regular polygons), we usually
consider only the lowest few.

An orbit $G(\lambda)$ of a Coxeter group $G$ is the set of points in
$\R^n$ generated by $G$ from a single seed point $\lambda\in\R^n$.
$G$-orbits are not common objects in the literature, nor is their
multiplication, which can be viewed in parallel to the
multiplication of $G$-invariant polynomials $P(\lambda;x)$
introduced in subsection 6.1 (for more about the polynomials see
\cite{D} and the references therein).\footnote{Polynomials in 6.1 are the simplest $W$-invariant ones. We are not concerned about any other of  their properties.} Indeed, the set of exponents
of all the monomials in $P(\lambda;x)$ is the set of points of the
orbit $G(\lambda)$. 

In this paper, we have adopted a point of view according to which the orbits $G(\lambda)$, being simpler than the polynomials $P(\lambda;x)$ or the weight systems of representations, are the primary objects of study.

The relation between the orbits of $W$ and the weight systems of
finite dimensional irreducible representations of semisimple Lie
groups/algebras over $\C$, can be understood as follows. The character
of a particular representation involves summation over the weight
system of the representation, i.e. over several $W$-orbits. As for which
orbits appear in a particular representation, this is a well known
question about multiplicities of dominant weights. There is a
laborious but rather fast computer algorithm for calculating the
multiplicities. Extensive tables of multiplicities can be found in
\cite{BMP}; see also the references therein. Thus one is justified in
assuming that the relation between a representation and a particular
$W$-orbit is known in all cases of interest.

Numerical characteristics, such as congruence classes, indices of
various degrees, and anomaly numbers, introduced here for
$W$-orbits, mirror similar properties of weight systems from
representation theory, which are often used in applications (for
example \cite{GPS,ST1,ST2,McP,S}).

In this paper, we introduce operations on $W$-orbits that are well
known for weight systems of representations: (i) The product of
$W$-orbits (of the same group) and its decomposition into the sum of
$W$-orbits; (ii) The decomposition of the $k$-th power of a
$W$-orbit symmetrized by the group of permutations of $k$ elements.
New is the introduction of such operations for the orbits of
non-crystallographic Coxeter groups. We intend to describe
reductions of $G$-orbits to orbits of a subgroup $G'\subset G$
 in a separate paper \cite{HLP}. Again, the
involvement of non-crystallographic groups makes the reduction
problem rather unusual. Corresponding applications deserve to be explored.

The decomposition of products of orbits of Coxeter groups, as
introduced here, is the core of other decomposition problems in
mathematics, such as the decomposition of direct products of
representations of semisimple Lie groups, the decomposition of
products of certain special functions \cite{P} and the decomposition
of products of $G$-invariant polynomials of several variables
\cite{D}. The last two problems are completely solved in terms
of orbit decompositions. The first problem requires that the
multiplicities of dominant weights in weight systems of
representations \cite{BMP} be known.

We view the $G$-orbits from a perspective uncommon in the
literature. Namely, the points of a $G$-orbit are taken to be
vertices of an $n$-dimensional $G$-invariant polytope centered at
origin, $n$ being the number of elementary reflections generating
$G$ (at the same time it is the rank of the corresponding semisimple
Lie group). The multiplication of two such polytopes/orbits, say
$P_1$ and $P_2$, is the set of points/vertices obtained by adding to
every point of $P_1$ every point of $P_2$. The resulting set of
points is again $G$-invariant and thus it is a union (we say `sum')
of several $G$-orbits (we say `$G$-polytopes'). Thus we have a ring
of $G$-polytopes with positive integer coefficients. We recall and
illustrate a general method of description of $n$-dimensional
reflection-generated polytopes \cite{MP2,CKPS}.

The core of our geometric interpretation of orbits as polytopes is
in the paragraph following equation \eqref{decomposition}. A product
of orbits is a union of concentric orbits. Geometrically this can be
seen as an  `onion'-like layered structure of orbits of different
radii. Unlike in representation theory, where orbit points are
always points of the corresponding weight lattices, in our case the
seed point of an orbit can be anywhere in $\R^n$. In particular, a
suitable choice of the seed points of the orbits, which are being
multiplied, can bring some of the layers of the `onion' structure as
close or as far apart as desired. Two examples are given in the last
section (see \eqref{close} and \eqref{thiny}).

\section{Reflections generating finite Coxeter groups}

Let $\a$ and $x$ be vectors in $\R^n$. We denote by $r_\a$ the
reflection in the $(n-1)$-dimensional `mirror' orthogonal to $\a$
and passing through the origin. For any $x\in \R^n$, we have
\begin{equation}\label{reflection}
 r_\a x=x-\frac{2\l x,\a\r}{\l\a,\a\r}\a\,.
\end{equation}
Here $\l a,b\r$ denotes scalar product in $\R^n$. In particular, we
have $r_\a 0=0$ and $r_\a\a=-\a$ so that $r^2_\a=1$.

A Coxeter group $G$ is by definition generated by several
reflections in mirrors that have the origin as their common point.
Various Coxeter groups are thus specified by the set $\Pi(\a)$ of
vectors $\a$, orthogonal to the mirrors and called the simple roots of $G$. Consequently, $G$ is given once the relative angles between elements of  $\Pi(\a)$ are given.

A standard presentation of $G$, generated by $n$ reflections, amounts
to the following relations
$$
r^2_k=1\,,\qquad (r_ir_j)^{m_{ij}}=1\,,\qquad k,i,j\in\{1,\dots,n\}\,,
$$
where we have simplified the notation by setting $r_{\a_k}=r_k$, and
where $m_{ij}$ are the lowest possible positive integers. The matrix
$(m_{ij})$ specifies the group. The angles between the mirrors of
reflections $r_i$ and $r_j$ are determined from the values of the exponents
$m_{ij}$. Indeed, for $m_{ij}=p$, the angle is $\pi/p$, while the angle between $\a_i$ and $\a_j$ is $\pi-\pi/p$.

The classification of finite reflection (Coxeter) groups was
accomplished in the first half of the 20th century.

\subsection{$\textbf{\emph{n}=1}$}\

There is just one group of order 2. Its two elements are $1$ and
$r$. We denote this group by $A_1$. Acting on a point $a$ of the
real line, the group $A_1$ generates its orbit of two points, $a$
and $ra=-a$, except if $a=0$. Then the orbit consists of just one
point, namely the origin.

\subsection{$\textbf{\emph{n}=2}$}\

There are infinitely many Coxeter groups in $\R^2$, one for each
$m_{12}=2,3,4,\dots$. Their orders are $2m_{12}$. In physics
literature, these are the dihedral groups.

Note that for $m_{12}=2$, the group is a product of two groups from
$n=1$. The reflection mirrors are orthogonal.

Our notation for the lowest five groups, generated by two
reflections, and their orders, is as follows:
\begin{alignat*}{2}
&m_{12}=2\quad :\quad A_1\times A_1\,, \quad &&\ 4 \\
&m_{12}=3\quad :\quad  A_2\,,          \quad &&\ 6 \\
&m_{12}=4\quad :\quad C_2\,,           \quad &&\ 8 \\
&m_{12}=5\quad :\quad H_2\,,           \quad &&10 \\
&m_{12}=6\quad :\quad G_2\,,           \quad &&12\,.
\end{alignat*}

\subsection{General case: Coxeter and Dynkin diagrams}\

A convenient general way to provide a specific set $\Pi(\a)$ is to draw a
graph where vertices are traditionally shown as small circles, one
for each $\a\in\Pi$, and where edges indicate absence of
orthogonality between two vertices linked by an edge.

A diagram consisting of several disconnected components means that
the group is a product of several pairwise commuting subgroups. Thus
it is often sufficient to consider only the groups with connected
diagrams.

In this paper, a Coxeter diagram is a graph providing only relative
angles between simple roots while ignoring their lengths. This is done by
writing $m_{ij}$ over the edges of the diagram. By convention, the
most frequently occurring value, $m_{ij}=3$, is not shown in the
diagrams. When $m_{ij}=2$, the edge is not drawn, i.e. the nodes
numbered $i$ and $j$ are not directly connected.

Consider the examples of Coxeter diagrams of all finite
non-crystallographic Coxeter groups with connected diagrams. Note
that we simply write $H_2$ when $m=5$.

\parbox{.6\linewidth}{\setlength{\unitlength}{2pt}
\def\kr{\circle{4}}
\thicklines
\begin{picture}(180,20)

\put(0,10){\makebox(0,0){${H_4}$}}
\put(10,10){\kr}
\put(20,10){\kr}
\put(30,10){\kr}
\put(40,10){\kr}
\put(12,10){\line(1,0){6}}
\put(22,10){\line(1,0){6}}
\put(32,10){\line(1,0){6}}
\put(34,11){$5$}

\put(60,10){\makebox(0,0){${H_3}$}}
\put(70,10){\kr}
\put(80,10){\kr}
\put(90,10){\kr}
\put(72,10){\line(1,0){6}}
\put(82,10){\line(1,0){6}}
\put(84,11){$5$}

\put(112,10){\makebox(0,0){${H_2(m)}$}}
\put(125,10){\kr}
\put(135,10){\kr}
\put(127,10){\line(1,0){6}}
\put(128,11){$m$}
\put(142,9){$m=5,7,8,9,\dots$}

\end{picture}}

A Dynkin diagram is a graph providing, in addition to the relative
angles, the relative lengths of the vectors from $\Pi(\a)$. Dynkin
diagrams are used for the crystallographic Coxeter groups,
frequently called the Weyl groups. There are four infinite series of
classical groups and five isolated cases of exceptional simple Lie
groups. Here is a complete list of Dynkin diagrams of such groups
(with connected diagrams):

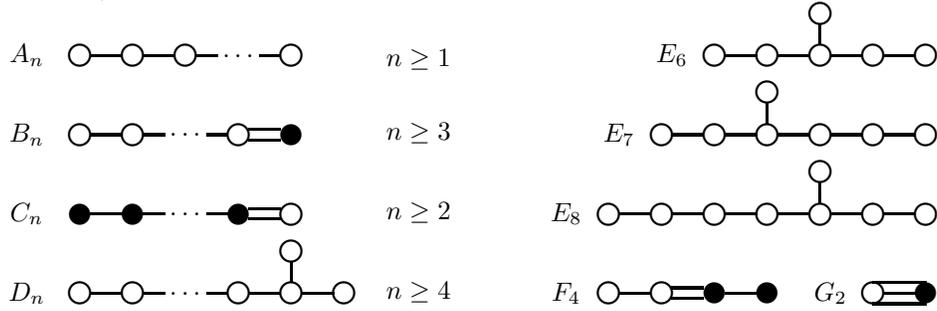
\parbox{.6\linewidth}{\setlength{\unitlength}{2pt}
\def\kr{\circle{4}}
\def\cr{\circle*{4}}
\thicklines
\begin{picture}(180,65)

\put( 0,55){\makebox(0,0){${A_n}$}}
\put(10,55){\kr}
\put(20,55){\kr}
\put(30,55){\kr}
\put(37,54.7){$\dots$}
\put(50,55){\kr}
\put(12,55){\line(1,0){6}}
\put(22,55){\line(1,0){6}}
\put(32,55){\line(1,0){4}}
\put(44,55){\line(1,0){4}}
\put(68,53){$n\geq 1$}

\put( 0,40){\makebox(0,0){${B_n}$}}
\put(10,40){\kr}
\put(20,40){\kr}
\put(27,39.7){$\dots$}
\put(40,40){\kr}
\put(50,40){\cr}
\put(12,40){\line(1,0){6}}
\put(22,40){\line(1,0){4}}
\put(34,40){\line(1,0){4}}
\put(42,41){\line(1,0){6}}
\put(42,39){\line(1,0){6}}
\put(68,39){$n\geq 3$}

\put( 0,25){\makebox(0,0){${C_n}$}}
\put(10,25){\cr}
\put(20,25){\cr}
\put(40,25){\cr}
\put(50,25){\kr}
\put(12,25){\line(1,0){6}}
\put(22,25){\line(1,0){4}}
\put(27,24.7){$\dots$}
\put(34,25){\line(1,0){4}}
\put(42,26){\line(1,0){6}}
\put(42,24){\line(1,0){6}}
\put(68,24){$n\geq 2$}

\put( 0,10){\makebox(0,0){${D_n}$}}
\put(10,10){\kr}
\put(20,10){\kr}
\put(40,10){\kr}
\put(50,10){\kr}
\put(60,10){\kr}
\put(50,18){\kr}
\put(12,10){\line(1,0){6}}
\put(22,10){\line(1,0){4}}
\put(27,9.7){$\dots$}
\put(34,10){\line(1,0){4}}
\put(42,10){\line(1,0){6}}
\put(52,10){\line(1,0){6}}
\put(50,12){\line(0,1){4}}
\put(68,9){$n\geq 4$}

\put(122,55){\makebox(0,0){${E_6}$}}
\put(130,55){\kr}
\put(140,55){\kr}
\put(150,55){\kr}
\put(160,55){\kr}
\put(170,55){\kr}
\put(150,63){\kr}
\put(132,55){\line(1,0){6}}
\put(142,55){\line(1,0){6}}
\put(152,55){\line(1,0){6}}
\put(162,55){\line(1,0){6}}
\put(150,57){\line(0,1){4}}

\put(112,40){\makebox(0,0){${E_7}$}}
\put(120,40){\kr}
\put(130,40){\kr}
\put(140,40){\kr}
\put(150,40){\kr}
\put(160,40){\kr}
\put(170,40){\kr}
\put(140,48){\kr}
\put(122,40){\line(1,0){6}}
\put(132,40){\line(1,0){6}}
\put(142,40){\line(1,0){6}}
\put(152,40){\line(1,0){6}}
\put(162,40){\line(1,0){6}}
\put(140,42){\line(0,1){4}}

\put(102,25){\makebox(0,0){${E_8}$}}
\put(110,25){\kr}
\put(120,25){\kr}
\put(130,25){\kr}
\put(140,25){\kr}
\put(150,25){\kr}
\put(160,25){\kr}
\put(170,25){\kr}
\put(150,33){\kr}
\put(112,25){\line(1,0){6}}
\put(122,25){\line(1,0){6}}
\put(132,25){\line(1,0){6}}
\put(142,25){\line(1,0){6}}
\put(152,25){\line(1,0){6}}
\put(162,25){\line(1,0){6}}
\put(150,27){\line(0,1){4}}

\put(102,10){\makebox(0,0){${F_4}$}}
\put(110,10){\kr}
\put(120,10){\kr}
\put(130,10){\cr}
\put(140,10){\cr}
\put(112,10){\line(1,0){6}}
\put(122,11){\line(1,0){6}}
\put(122,9){\line(1,0){6}}
\put(132,10){\line(1,0){6}}

\put(152,10){\makebox(0,0){${G_2}$}}
\put(160,10){\kr}
\put(170,10){\cr}
\put(162,10){\line(1,0){6}}
\put(160,12){\line(1,0){10}}
\put(160,8){\line(1,0){10}}

\end{picture}}

The names of the groups, as is traditional in Lie theory, are shown
on the left of each diagram. Open (black) circles indicate longer
(shorter) roots. The ratio of their square lengths is
$\l\a_l,\a_l\r:\l\a_s,\a_s\r=2:1$ in all cases except for $G_2$
where the ratio is $3:1$. Moreover, we adopt the usual convention
that $\l\a_l,\a_l\r=2$. A single, double, and triple line indicates
respectively the angle $2\pi/3$, $3\pi/4$, and $5\pi/6$ between the
roots, or equivalently, the angles $\pi/3$, $\pi/4$, and $\pi/6$
between the reflection mirrors. The absence of a direct link between
two nodes implies that the corresponding simple roots, as well as
the mirrors, are orthogonal. Note that the relative angles of the
mirrors of $B_n$ and $C_n$ coincide. Hence their $W$-groups are
isomorphic. Their simple roots differ by length.

We adopt the Dynkin numbering of nodes. The numbering
proceeds from left to right $1,2,\dots$ In case of $D_n$ and
$E_6$, $E_7$, $E_8$, the node above the main line has the highest
number, respectively $n$, 6, 7, 8.

Orders of the finite Coxeter groups are provided in Table~1 for
groups with connected diagrams. When a diagram has several
disconnected components, the order is the product of orders
corresponding to each subdiagram.


\begin{table}
\begin{center}
\begin{tabular}{|c|c|c|c|c|c|}
\hline
$A_n\ (n\geq1)$&$B_n\ (n\geq3)$&$C_n\ (n\geq2)$&$D_n\ (n\geq4)$&$E_6$&$E_7$\\
\hline
$(n+1)!$&$2^nn!$&$2^nn!$&$2^{n-1}n!$&$2^73^45$&$2^{10}3^457$\\
\hline
\hline
$E_8$&$F_4$&$G_2$&$H_2(m)$&$H_3$&$H_4$\\
\hline
$2^{14}3^55^27$&$2^73^2$&$12$&$2m$&$120$&$120^2$\\
\hline
\end{tabular}
\medskip

\caption{Orders of the finite Coxeter groups.}
\label{tab1}
\end{center}
\end{table}

\section{Root and weight lattices}

Information essentially equivalent to that provided by the Coxeter
and Dynkin diagrams is also given in terms of $n\times n$ matrices
$C$, called the Cartan matrices. Relative angles and lengths of
simple roots can be used to form the Cartan matrix for each group.
Its matrix elements are calculated as
\begin{equation}\label{cartan}
C=(C_{jk})=\left(\frac{2\l\a_j,\a_k\r}{\l\a_k,\a_k\r}\right)\,,
\qquad j,k\in\{1,2,\dots,n\}\,.
\end{equation}
Cartan matrices and their inverses are given in many places, e.g.
\cite{H,BMP}.

The Cartan matrices can be defined for any finite
Coxeter group by using formula~\eqref{cartan}. For
non-crystallographic groups the matrices are
$$
C(H_2)=\left( \begin{smallmatrix} 2 &-\tau \\ -\tau & 2
\end{smallmatrix}\right)\,,\quad C(H_3)=\left( \begin{smallmatrix}
2&-1&0\\ -1&2&-\tau\\ 0&-\tau&2\end{smallmatrix}\right)\,,\quad
C(H_4)=\left( \begin{smallmatrix} 2&-1&0&0\\ -1&2&-1&0\\
0&-1&2&-\tau\\ 0&0&-\tau&2\end{smallmatrix}\right)\,,
$$
where $\tau$ is the larger of the solutions of the algebraic
equation $x^2 = x + 1$, i.e. $\tau=\tfrac12(1+\sqrt5)$.

In addition to the basis of simple roots ($\a$-basis), it is useful
to introduce the basis of fundamental weights ($\o$-basis).
Subsequently, most of our computations will be performed in the
$\o$-basis.
$$
\a=C\o\,,\qquad \o=C^{-1}\a\,.
$$
Note the important relation:
\begin{equation}\label{ortohog}
\l\a_k,\o_j\r=\delta_{jk}\frac{\l\a_k,\a_k\r}2\,,\qquad
    j,k\in\{1,2,\dots,n\}\,.
\end{equation}

Illustrations showing the $\a$- and $\o$-bases of $A_2$, $C_2$, and
$G_2$ are given in Figure~1 of \cite{MP2}.

The root lattice $Q$ and the weight lattice $P$ of $G$ are formed by
all integer linear combinations of simple roots, respectively
fundamental weights, of $G$,
\begin{equation}\label{lattices}
Q=\Z\a_1+\cdots+\Z\a_n\,,\qquad P=\Z\o_1+\cdots+\Z\o_n\,.
\end{equation}
Here $\Z$ stands for any integer. For the groups that have simple roots of two different lengths, one may define the root lattice of $n$ linearly independent short roots, which cannot all be simple. In general, $Q\subseteq
P$, with $Q = P$ only for $E_8$, $F_4$, and $G_2$.

If $G$ is one of the non-crystallographic Coxeter groups, the
lattices $Q$ and $P$ are dense everywhere.

Since $\a$- and $\o$-bases are not orthogonal and not normalized, it
is sometimes useful to work with orthonormal bases. For
crystallographic groups, they are found in many places, for example
\cite{Bourbaki,BMP}. For non-crystallographic groups, $H_2$, $H_3$
and $H_4$; see \cite{CMP,MP93}.

\section{The orbits of Coxeter groups}

\subsection{Computing points of an orbit}\

Given the reflections $r_\alpha$, $\alpha\in\Pi(\alpha)$, of a
Coxeter group $G$, and a seed point $\lambda\in\R^n$, the points of
the orbit $G(\lambda)$ are given by the set of distinct points
generated by repeated application of the reflections $r_\alpha$ to
$\lambda$. All points of an orbit are equidistant from the origin.
The radius of an orbit is the distance of (any) point of the orbit
from the origin.

There are practically important considerations which make it almost
imperative that the computation of the points of any orbit of $G$ be
carried out in the $\o$-basis, as follows:

\begin{itemize}
\item
Every orbit contains precisely one point with nonnegative coordinates
in the $\o$-basis. We specify the orbit by that point, calling it the
dominant point of the orbit.\\
\item
Given a dominant point $\lambda$ of the group $G$ in the $\o$-basis, one
readily finds the size of the orbit $G(\lambda)$, i.e. the number of
points in the orbit, using the order $|G|$ of the Coxeter group and
the order of the stabilizer of $\lambda$ in $G$:
\begin{equation}\label{orbitsize}
|G(\lambda)|=\frac{|G|}{|\operatorname{Stab}_G(\lambda)|}
\end{equation}
Here $\operatorname{Stab}_G(\lambda)$ is a Coxeter subgroup of $G$.
To find it, one needs to attach the $\o$-coordinates of $\lambda$ to
the corresponding nodes of the diagram of $G$. The subdiagram
carrying the coordinates 0 is the diagram of  $\operatorname{Stab}_G(\lambda)$.\\
\item
Due to \eqref{ortohog}, the reflections \eqref{reflection} are
particularly simple when applied to $\o$'s:
\begin{equation}\label{michelle}
r_k\o_j=\o_j-\frac{2\l\a_k,\o_j\r}{\l\a_k,\a_k\r}\a_k
       =\o_j-\delta_{jk}\a_k\,.
\end{equation}\\
\item
Starting from the dominant point of an orbit, it suffices to apply,
during the computation of the orbit points, only reflections corresponding
to positive coordinates of any given weight. All points of the orbit are found in this way.\\

\end{itemize}

\subsection{Orbits of $A_2$, $C_2$, $G_2$, and $H_2$}\

We give some examples of orbits. Let $a,b > 0$.
\begin{align*}
A_2\ :\quad   G((a,0)) = & \{(a,0),\, (-a,a),\, (0,-a)\}\,,\\
                G((0,b)) = & \{(0,b),\, (b,-b),\, (-b,0)\}\,,\\
                G((a,b)) = & \{(a,b),\, (-a,a+b),\, (b,-a-b),\,(a+b,-b)\,,\\
                        \notag & (-a-b,a),\, (-b,-a)\}\,.\\
\intertext{In particular, the orbit $G((1,1))=\{(1,1),\, (-1,2),\,
(1,-2),\, (2,-1),\, (-2,1),\, (-1,-1)\}$ consists of the vertices of
a regular hexagon of radius $\sqrt2$. It is the root system of $A_2$.}
C_2\ :\quad   G((a,0)) = & \{\pm(a,0),\, \pm(-a,a)\}\,,\\
                G((0,b)) = & \{\pm(0,b),\, \pm(2b,-b)\}\,,\\
                G((a,b)) = & \{\pm(a,b),\, \pm(-a,a+b),\, \pm(a+2b,-a-b),\, \\
                        \notag & \pm(-a-2b,b)\}.
\intertext{In particular, the orbits $G((2,0))$ and $G((0,1))$ of
radii $\sqrt2$ and $1$ are respectively the vertices and midpoints
of the sides of a square. Together the two orbits form the root system of $C_2$.}
G_2\ :\quad G((a,0)) =&\{\pm(a,0),\,\pm(-a,3a),\,\pm(2a,-3a)\}\,,\\
            G((0,b)) =&\{\pm(0,b),\,\pm(b,-b),\,\pm(-b,2b)\}\,,\\
            G((a,b)) =&\{\pm(a,b),\,\pm(-a,3a+b),\,\pm(2a+b,-3a-b),\,\\
            \notag&\pm(-2a-b,3a+2b),\,\pm(a+b,-3a-2b)\,,\\
            \notag&\pm(-a-b,b)\}\,.
\intertext{In particular, the orbits $G((1,0))$ and $G((0,1))$ are the
vertices of regular hexagons of radii $\sqrt2$ and $2/\sqrt3$,
rotated relatively by $30^\circ$, i.e. they form a hexagonal  star. Together the two orbits form the root system of $G_2$.
The points of  $G((a,a\sqrt3/\sqrt2))$, $a>0$, are the vertices of a
regular dodecahedron of radius $\sqrt2a$.}
H_2\ :\quad G((a,0))
=&\{(a,0),\, (-a,a\tau),\, (a\tau,-a\tau),\,
                        (-a\tau,a),\, (0,-a)\}\,,\label{(a,0)H2}\\
             G((0,b)) =&\{(0,b),\, (b\tau,-b),\, (-b\tau,b\tau),\,
                        (b,-b\tau),\, (-b,0)\}\,,\\
             G((a,b)) =&\{(a,b),\, (-a,b+a\tau),\, (a\tau +b\tau,-b-a\tau),\, \\
                        \notag&  (-a\tau-b\tau,a+b\tau),\, (b,-a-b\tau),\, (a+b\tau,-b),\, \\
                        \notag& (-a-b\tau,a\tau+b\tau),\, (b+a\tau,-a\tau-b\tau),\,\\
                        \notag& (-b-a\tau,a),\, (-b,-a)\}\,.
\end{align*}
In particular, the orbits $G((a,0))$ and $G((0,b))$ are the vertices of
regular pentagons of radii $a\sqrt2$ and $b\sqrt2$, rotated
relatively by $36^\circ$. The orbit $G((a,a))$ forms a regular
decahedron. The orbit $G((\tau,\tau))$ consists of the roots of $H_2$.

An orbit of $A_2$ or $H_2$ contains, with every point $(p,q)$ also
the point $(-q,-p)$. Note that in the examples of this subsection
the constants $a$ and $b$ do not need to be integers. All one
requires is that they are positive. Effects of special choices of
these constants are exemplified in \eqref{close} and \eqref{thiny}
below.

\begin{table}
\begin{center}
\begin{tabular}{|c|c||c|c|c||c|c||c|c|c|}
\hline
&&$A_2$&$C_2$&$G_2$&$H_2$&$H_2(7)$  &$1$&$2$&$3$\\
\hline\hline
1&$\bullet\bullet$&$6$&$8$&$12$&$10$&$14$  &$\checkmark$&$$&$$\\
\hline
2&$\bullet\circ$&$3$&$4$&$6$&$5$&$7$  &$$&$\checkmark$&$$\\
\hline
3&$\circ\bullet$&$3$&$4$&$6$&$5$&$7$  &$$&$$&$\checkmark$\\
\hline \hline\hline
4&$\star\bullet$&$3$&$4$&$6$&$5$&$7$  &$\checkmark$&$\checkmark$&$$\\
\hline
5&$\bullet\star $&$3$&$4$&$6$&$5$&$7$  &$\checkmark$&$$&$\checkmark$\\
\hline
\end{tabular}
\medskip

\caption{Number of faces of $2D$ polytopes with Coxeter group
symmetry. The first three rows specify representatives of $G$-orbits
of $2D$ polytopes. A black (open) dot in the second column stands
for a positive (zero) coordinate in the $\o$-basis of the dominant
point representing the orbit of vertices. The number of vertices is
listed in the subsequent five columns. Rows 4 and 5 refer to the
edges of the polytopes. A star in the second column indicates the
reflection generating the symmetry group of the edge. The number of
edges is shown for each group in subsequent columns. Check marks in
one of the last three columns indicate the faces which belong to the
polytope described in that column.} \label{tab2}
\end{center}
\end{table}

\eject
\subsection{Orbits of $A_3$, $B_3$, $C_3$, and $H_3$}\

We give some examples of orbits. Let $a,b > 0$.

\begin{align*}
A_3\ :\quad G((a,0,0))= &\{(a,0,0),\, (-a,a,0),\, (0,-a,a),\,
                        (0,0,-a)\}\,,\\
            G((0,b,0))= &\{\pm(0,b,0),\,\pm(b,-b,b),\,\pm(-b,0,b)\}\,,\\
            G((1,1,0))= &\{(1,1,0),\, (-1,2,0),\, (2,-1,1),\, (1,-2,2),\,
                        (-2,1,1),\,\\
                        \notag & (2,0,-1),\, (-1,-1,2),\, (1,0,-2),\,
                        (-2,2,-1),\,\\
                        \notag & (-1,1,-2),\, (0,-2,1),\, (0,-1,-1)\}\,.\\
\intertext{}
B_3\ :\quad G((a,0,0))= &\{\pm(a,0,0),\, \pm(-a,a,0),\,
             \pm(0,-a,2a)\}\,,\\
            G((0,b,0))= &\{\pm(0,b,0),\, \pm(b,-b,2b),\,
            \pm(-b,0,2b),\, \pm(b,b,-2b),\, \\
            \notag & \pm(-b,2b,-2b),\, \pm(2b,-b,0) \}\,,\\
            G((0,0,c))= &\{\pm(0,0,c),\, \pm(0,c,-c),\,
            \pm(c,-c,c),\, \pm(c,0,-c)\}\,. \\
\intertext{}
C_3\ :\quad G((a,0,0))= &\{\pm(a,0,0),\, \pm(-a,a,0),\, \pm(0,-a,a)\}\,,\\
            G((0,b,0))= &\{\pm(0,b,0),\, \pm(b,-b,b),\, \pm(-b,0,b),\,
                        \pm(b,b,-b),\,\\
                        \notag & \pm(-b,2b,-b),\, \pm(2b,-b,0)\}\,,\\
            G((0,0,c))= &\{\pm(0,0,c),\, \pm(0,2c,-c),\,
                        \pm(2c,-2c,c),\, \pm(2c,0,-c)\}\,. \\
\intertext{}
H_3\ :\quad G((a,0,0))= &\{\pm(a,0,0),\, \pm(-a,a,0),\,\pm(0,-a,a\tau),\,
                        \pm(0,a\tau,-a\tau),\\
                        \notag & \pm(a\tau,-a\tau,a),\,\pm(-a\tau,0,a)\}\,,\\
            G((0,0,c))= &\{\pm(0,0,c),\, \pm(0,\tau c,-c),
                         \pm(\tau c,-\tau c,\tau^2 c),\pm(-\tau c,0,\tau^2 c),\\
                        \notag & \pm(a\tau,-a\tau,a),\,\pm(-a\tau,0,a)\}\,.
\end{align*}

\begin{table}
\begin{center}
\begin{tabular}{|c|c||c|c|c||c||c|c|c|c|c|c|c|}
\hline
&Diagram&$A_3$&$B_3$&$C_3$&$H_3$       &1&2&3&4&5&6&7\\
\hline\hline
1&$\bullet\bullet\bullet$&24&48&48&120 &$\checkmark$&&&&&&\\
\hline
2&$\bullet\bullet\circ$&12&24&24&60    &&$\checkmark$&&&&&\\
\hline
3&$\bullet\circ\bullet$&12&24&24&60    &&&$\checkmark$&&&&\\
\hline
4&$\circ\bullet\bullet$&12&24&24&60  &&&&$\checkmark$&&&\\
\hline
5&$\bullet \circ \circ $&4&6&6&12  &&&&&$\checkmark$&&\\
\hline
6&$\circ\bullet \circ $&6&12&12&30   &&&&&&$\checkmark$&\\
\hline
7&$\circ \circ\bullet$&4&8&8&20  &&&&&&&$\checkmark$\\
\hline \hline \hline
8&$\star\bullet\bullet$&12&24&24&60   &$\checkmark$&&$\checkmark$&&&&\\
\hline
9&$\bullet\star\bullet$&12&24&24&60 &$\checkmark$&$\checkmark$&&$\checkmark$&&$\checkmark$&\\
\hline
10&$\bullet\bullet\star$&12&24&24&60  &$\checkmark$&&$\checkmark$&&&&\\
\hline
11&$\star\bullet\circ$&6&12&12&30     &&$\checkmark$&&&$\checkmark$&&\\
\hline
12&$\circ\bullet\star$&6&12&12&30     &&&&$\checkmark$&&&$\checkmark$\\
\hline
13&$\star\star\bullet$&4&8&8&20  &$\checkmark$&$\checkmark$&$\checkmark$&$\checkmark$&$\checkmark$&$\checkmark$&\\
\hline
14&$\star\bullet\star$&6&12&12&30     &$\checkmark$&&$\checkmark$&&&&\\
\hline
15&$\bullet\star\star$&4&6&6&12 &$\checkmark$&$\checkmark$&$\checkmark$&$\checkmark$&&$\checkmark$&$\checkmark$\\
\hline
\end{tabular}
\medskip

\caption{The $3D$ polytopes generated by a Coxeter group from a single
seed point, and the number of their faces of dimension 0, 1 and 2.
Decorated diagrams, rows 1 to 7, specify the polytopes. The dimension
of a face equals the number of stars in the diagram. See 5.1 for additional explanations.} \label{tab3}
\end{center}
\end{table}

\section{Orbits as polytopes}

In this section, we recall an efficient method \cite{CKPS} of
description for reflection-generated polytopes in any dimension.

The idea of the method consists in the following. Suppose we have an
orbit $G(\lambda)$. Consider its points as vertices (faces of
dimension 0) of the polytope also denoted $G(\lambda)$ in $\R^n$.
Then for any face $f$ of dimension $0\leq d\leq n-1$, we identify
its stabilizer $\operatorname{Stab}_{G(\lambda)}(f)$ in $G$, which
is a product of two Coxeter subgroups of $G$:
$$
\operatorname{Stab}_{G(\lambda)}(f)= G_1(f)\times G_2(D)
$$
where $G_1(f)$ is the symmetry group of the face, and $G_2(D)$ stabilizes
$f$ pointwise, i.e. does not move it at all.

Our method consists in recursive decorations of the diagram of $G$,
providing at each stage the subdiagrams of $G_1(f)=G(\star)$ and
$G_2(D)=G(\circ)$ for faces of one type. The decoration of the nodes
of the diagram indicates to which $G(\star)$ or $G(\circ)$ subgroups
of the stabilizer the corresponding reflections belong. For further
details, see \cite{CKPS}. A much wider application of this method is
described in \cite{MP2,MP3,MP4}, including its exploitation in non
Euclidean spaces.

We start with an extreme decoration of the diagram. It is
equivalent to stating which coordinates of the dominant weight are
positive relative to the $\o$-basis. The nodes are drawn as either
open or black circles, i.e. zero or positive coordinates
respectively.

Every possible extreme decoration fixes a polytope. There are only
two rules for recursive decoration of the diagrams, starting from
one of the extreme ones: (i) A single black circle is
replaced by a star; (ii) open circles, that become
adjacent to a star by diagram connectivity, are changed to black
ones.

Tables~2 and 3 show the results of the application of the decoration
rules for polytopes in $2D$ and $3D$ for all groups with connected
diagrams. All polytopes for $A_4$, $B_4$, $C_4$, $D_4$, and
$H_4$ are described in Tables~3 and 4 of \cite{CKPS}.

\subsection{Explanation of the Tables}\

A description of Table~2 is given in its caption.

Consider Table~3. The second column contains short-hand notation for
several diagrams at once. We call them decorated diagrams. No links
between nodes of a diagram are drawn because they would need to be
different for each group in subsequent columns. The nodes do not reveal the relative lengths
of roots, their decoration indicates to which of the pertinent
subgroup of the stabilizer of $G$ such a reflection belongs. Thus the diagrams of the
second column of the Table apply to  $A_3$, $B_3$, $C_3$ and $H_3$
at the same time.

Each line of the Table describes one of $G$-orbits of identical
faces. The dimension of the face equals to the number of stars in its
decorated diagram. Numerical entries in a row give the number of
faces for polytopes of symmetry groups $A_3$, $B_3$, $C_3$ and
$H_3$, shown in the header of the columns. The top seven rows show
the starting decorations fixing the polytopes, and also the number
of 0-faces (vertices) of the polytopes of each group. The check
marks in one of the last seven columns indicate the faces belonging
to the same polytope.

 \begin{example}\

As an example of how to decipher properties of polytope faces, consider
rows number 5 and 2. The diagram in row 5 conveys the fact that
$\lambda=a\o_1$ with $a>0$. The exact value of $a$ affects only the
size of the polytope, not its shape. The stabilizer of $\lambda$ is
given by the subdiagram of open circles, i.e. $r_2$ and $r_3$
generate its stabilizer. For $A_3$ the subdiagram is of type $A_2$,
while for $B_3$ and $C_3$ it is of type $C_2$, and for $H_3$ it is
of type $H_2$. Hence in row 5 the entries give the number of
vertices as $24/6$, $48/8$, $48/8$, $120/10$ respectively.

The check mark in column 5 and row 5 indicates that faces belonging
to our polytope are indicated by other check marks in column 5,
namely in rows 11 and 13. The diagram of row 11 has just one star,
hence the face is 1-dimensional (an edge). Its stabilizer (the
subdiagram of stars and open circles) is of type $A_1\times A_1$ for
all four cases. Hence the number of edges is $24/4$ for $A_3$,
$48/4$ for $B_3$ and $C_3$, and $120/4$ for $H_3$. The only type of
$2D$ face is given in row 13. The symmetry group of the face is generated
by $r_1$ and $r_2$. It is of type $A_2$ for all four cases. Thus
there are $24/6$ faces in $A_3$, $48/6$ in $B_3$ and $C_3$, and
$120/6$ in $H_3$ polytope.

Similarly, row 2 indicates that $\lambda=a\o_1+b\o_2$, $a,b>0$. It is stabilized by the group generated by $r_3$, which is of type $A_1$ for all four cases. Hence
the number of vertices equals half of the order of the corresponding
Coxeter group. There are two orbits of edges given in rows 9 and 11,
while the two orbits of $2D$ faces are given by the check marks in
rows 13 and 15.
 \end{example}

 \begin{example}\

The $2D$ faces can actually be constructed knowing their symmetry and
the seed point, say $(a,0,0)$. The diagram of the $2D$ face is $\star~\star~\bullet$, meaning that the symmetry group of the face is generated by $r_1$ and $r_2$. Moreover, it is of the same type ($A_2$) for all four groups. Then there are just three
distinct vertices of the $2D$ face:
$$
 (a,0,0),\quad r_1(a,0,0),\quad r_2r_1(a,0,0)\,.
 $$
The $2D$ face is formed from the seed point $(a,0,0)$ by application
of reflections $r_1$ and $r_2$.

The vertices of the $2D$ face are different triangles for each group,
because they are given in their respective $\o$-basis:
 \begin{align*}
 A_3\ :\quad &(a,0,0),\ (-a,a,0),\ (0,-a,a),\\
 B_3\ :\quad &(a,0,0),\ (-a,a,0),\ (0,-a,2a),\\
 C_3\ :\quad &(a,0,0),\ (-a,a,0),\ (0,-a,a),\\
 H_3\ :\quad &(a,0,0),\ (-a,a,0),\ (0,-a,a\tau)\,.
 \end{align*}
 \end{example}

 \begin{example}\

Let us consider row 2 in further detail. The starting point is
$\lambda=a\o_1+b\o_2$, where $a,b>0$. There are two orbits of edges given by their endpoints:
 $$
 ((a,b,0), r_1(a,b,0))\,,\qquad ((a,b,0), r_2(a,b,0))\,,
 $$
and two orbits of $2D$ faces. Consider just the $H_3$ case. The $2D$
face of row 13 has the symmetry group generated by $r_1,r_2$ ($A_2$
type). It is a hexagon:
\begin{align*}
 &(a,b,0),\quad (-a,a+b,0),\quad (a+b,-b,\tau b),\quad(b,-a-b,\tau(a+b)),\\
 &(-a-b,a,\tau b),\quad (-b,-a,\tau(a+b)).
\end{align*}
The $2D$ face of row 15 has its symmetry group generated by $r_2,r_3$
($H_2$ type). It is a pentagon:
\begin{align*}
 &(a,b,0),\quad (a+b,-b,\tau b),\quad (a+b,\tau b,-\tau b),\quad \\
 &(a+\tau^2 b, -\tau b,b),\quad (a+\tau^2 b,0,-b)\,.\\
\end{align*}
In particular, when $a=b$, the pentagon and the hexagon are both
regular. The polytope is then the familiar fullerene or `soccer
ball'.
  \end{example}

Further questions about the structure of polytopes can be answered
within our formalism: How many $2D$ faces meet in a vertex? Which
$2D$ faces meet in an edge? The higher the dimension, the more
questions like these can be asked and answered. For more information
on such questions and others (e.g. dual pairs of polytopes), we
refer to \cite{CKPS}.

\section{Decomposition of products of polytopes}

\subsection{Multiplication of $G$-invariant polynomials}\

The product of $G$-polytopes together with its decomposition, as
defined in 6.2 below, can be simply motivated by its correspondence
to the product of more familiar objects than orbits, namely
$G$-invariant polynomials, say $P(\lambda;x)$ and $P(\mu;x)$. Here
$\lambda$ and $\mu$ are dominant points of their orbits and $x$
stands for $n$ auxiliary independent variables $x_1,x_2,\dots,x_n$
whose nature is of no concern to us here. They can be thought of as,
for example, complex or real variables. We introduce them in order
to make sense of the definitions below.

Denote by $\lambda^{(i)}\in G(\lambda)$ the points of the orbit
$G(\lambda)$, and by $\mu^{(k)}\in G(\mu)$, where
\begin{equation}
\lambda^{(i)}=\sum_{p=1}^n a_p^{(i)}\o_p\,,\quad
\mu^{(k)}=\sum_{q=1}^n b_q^{(k)}\o_q\,,\qquad
1\leq i\leq|G(\lambda)|\,,\quad 1\leq k\leq|G(\mu)|\,.
\end{equation}
Here $|G(\lambda)|$ and $|G(\mu)|$ denote the number of points in
their orbits. Then we can introduce the polynomials:
\begin{align}\label{polynomial}
P(\lambda;x)&=\sum_{\lambda^{(i)}\in G(\lambda)}x^{\lambda^{(i)}}
           :=\sum_{i=1}^{|G(\lambda)|}
             x_1^{a_1^{(i)}}x_2^{a_2^{(i)}}\cdots x_n^{a_n^{(i)}},\\
P(\mu;x)&=\sum_{\mu^{(k)}\in G(\mu)}x^{\mu^{(k)}}
           :=\sum_{k=1}^{|G(\mu)|}
             x_1^{b_1^{(k)}}x_2^{b_2^{(k)}}\cdots x_n^{b_n^{(k)}}\,,
\end{align}
and their product,
\begin{equation}
P(\lambda;x)\otimes P(\mu;x)=\sum_{i=1}^{|G(\lambda)|}\sum_{k=1}^{|G(\mu)|}
       x_1^{a_1^{(i)}+b_1^{(k)}}x_2^{a_2^{(i)}+b_2^{(k)}}\cdots
       x_n^{a_n^{(i)}+b_n^ {(k)}}\,.
\end{equation}
The latter consists of the sum of $|G(\lambda)||G(\mu)|$ monomials
which can be decomposed into the sum of polynomials defined by one
$G$-orbit each.

Finally, consider an example: Let $G$ be the group $A_2$, and
$\lambda=(1,0)$ and $\mu=(0,1)$. Therefore
$P((1,0);x)=x_1+x_1^{-1}x_2+x_2^{-1}$ and
$P((0,1);x)=x_2+x_1x_2^{-1}+x_1^{-1}$. Their products decompose as
follows,
\begin{align*}
P((1,0);x)\otimes P((0,1);x)=&
\{x_1x_2+x_1^2x_2^{-1}+x_1^{-1}x_2^2+x_1^{-2}x_2+x_1x_2^{-2}
        +x_1^{-1}x_2^{-1}\}+3\\
                            =&P((1,1);x)+3P((0,0);x)\,,\\
P((1,0);x)\otimes P((1,0);x)=&
  \{x_1^2+x_1^{-2}x_2^2+x_2^{-2}\}+2\{x_2+x_1x_2^{-1}+x_1^{-1}\}\\
=& P((2,0);x)+2P((0,1);x)\,.
\end{align*}

\subsection{Products of $G$-orbits}\

Suppose we are given two orbits, say $G(\lambda)$ and $G(\mu)$, of
the same Coxeter group $G$. Let $\lambda^{(i)}$ and $\mu^{(k)}$ be the
points of $G(\lambda)$ and $G(\mu)$ respectively, numbered in some
way. We define the product of two orbits as
\begin{equation}\label{product}
G(\lambda)\otimes G(\mu)
  :=\bigcup_{\lambda^{(i)}\in G(\lambda),\ \mu^{(k)}\in G(\mu)}
  (\lambda^{(i)}+\mu^{(k)})\,.
\end{equation}
The left side is obviously $G$-invariant, therefore the right side is also
$G$-invariant. Hence it can be decomposed into a union of several
$G$-orbits. The highest and the lowest components of such a
decomposition are easily obtained:
\begin{equation}\label{decomposition}
G(\lambda)\otimes G(\mu)
     =G(\lambda+\mu)\cup\cdots\cup
      G(\lambda+\overline{\mu}).
\end{equation}
Here, $\lambda+\mu$ is the sum of the dominant points of the orbits
$G(\lambda)$ and $G(\mu)$. The symbol $\overline{\mu}$ stands for
the unique lowest point of $G(\mu)$ (all coordinates are
non-positive in the $\o$-basis). Frequently, it happens that
$\lambda+\overline{\mu}$ is not a dominant point, i.e. the highest
point in its orbit, but it still identifies the orbit uniquely. Note
also that $\lambda+\overline{\mu}$  and $\mu+\overline{\lambda}$
always belong to the same $G$-orbit. The lowest component often
appears more than once in the decomposition.

For a geometric interpretation of \eqref{decomposition}, recall that
all orbits in  \eqref{decomposition} are concentric, having the origin
as their common center, and that points of one orbit are equidistant
from the origin. In physics, the product on the left side of
\eqref{decomposition} can be thought of as a certain `interaction'
between two orbit-layers, resulting on the right side in an
`onion'-like structure of several concentric orbit-layers.

To simplify the notation in the following examples, we write just
$\lambda$ instead of $G(\lambda)$, so that $\lambda \otimes \mu$
means $G(\lambda)\otimes G(\mu)$.

\subsection{Two-dimensional examples}\
\begin{align*}
A_2\ :\quad (1,0)\otimes (0,1) &=(1,1)\cup 3(0,0)\,,\\
            (1,0)\otimes (1,1) &=(2,1)\cup 2(1,0)\cup 2(0,2)\,,\\
            (1,1)\otimes (1,1) &=(2,2)\cup 2(1,1)\cup 2(3,0)\cup 2(0,3)\cup 6(0,0)\,.
\intertext{}
C_2\ :\quad (1,0)\otimes (0,1) &=(1,1)\cup 2(1,0)\,,\\
            (1,0)\otimes (1,1) &=(2,1)\cup 2(2,0)\cup 2(0,2)\cup 2(0,1)\,,\\
            (1,1)\otimes (1,1) &=(2,2)\cup 2(2,1)\cup 2(4,0)\cup 2(2,0)\cup 2(0,3) \cup\\
                            \notag &\quad 2(0,1)\cup 8(0,0)\,.
\intertext{}
G_2\ :\quad (1,0)\otimes (0,1) &=(1,1)\cup 2(0,2)\cup 2(0,1)\,,\\
            (1,0)\otimes (1,1) &=(2,1)\cup (1,2)\cup (1,1)\cup 2(0,4)\cup 2(0,2)\cup 2(0,1)\,,\\
            (1,1)\otimes (1,1) &=(2,2)\cup 2(1,1)\cup 2(1,3)\cup 2(3,0)\cup 2(2,0) \cup\\
                            \notag &\quad 2(1,0)\cup2(0,5)\cup 2(0,4)\cup 2(0,1)\cup 12(0,0)\,.
\intertext{}
H_2\ :\quad (1,0)\otimes (0,1) &=(1,1)\cup (\tau-1,\tau-1)\cup 5(0,0)\,,\\
            (1,0)\otimes (1,1) &=(2,1)\cup (\tau,\tau-1)\cup (\tau-1,1)\cup2(1,0) \cup\\
                            \notag &\quad 2(0,\tau+1)\,,\\
            (1,1)\otimes (1,1) &=(2,2)\cup 2(\tau,\tau)\cup 2(\tau-1,\tau-1)\cup 2(2+\tau,0) \cup\\
                            \notag &\quad 2(2\tau-1,0)\cup 2(0,2+\tau)\cup 2(0,2\tau-1)\cup 10(0,0)\,.
\end{align*}

\subsection{Three-dimensional examples}\
\begin{align*}
A_3\ :\quad (1,0,0)\otimes (0,0,1) &=(1,0,1)\cup 4(0,0,0)\,,\\
            (1,0,1)\otimes (0,1,0) &=(1,1,1)\cup 3(2,0,0)\cup 4(0,1,0) \cup 3(0,0,2)\,,\\
            (1,1,0)\otimes (0,0,1) &=(1,1,1)\cup 3(2,0,0)\cup2(0,1,0)\,.
\intertext{}
B_3\ :\quad (1,0,0)\otimes (0,0,1) &=(1,0,1)\cup 3(0,0,1)\,,\\
            (1,0,1)\otimes (0,1,0) &=(1,1,1)\cup 2(2,0,1)\cup 3(1,0,1) \cup 2(0,1,1) \cup\\
                                \notag &\quad 3(0,0,3)\cup 6(0,0,1)\,,\\
            (1,1,0)\otimes (0,0,1) &=(1,1,1)\cup 2(2,0,1)\cup 2(1,0,1) \cup 2(0,1,1)\,.
\intertext{}
C_3\ :\quad (1,0,0)\otimes (0,0,1) &=(1,0,1)\cup 2(0,1,0)\,,\\
            (1,0,1)\otimes (0,1,0) &=(1,1,1)\cup 2(2,1,0)\cup 2(1,0,1) \cup 4(2,0,0) \cup\\
                                \notag &\quad 4(0,2,0)\cup 4(0,1,0)\cup 3(0,0,2)\,,\\
            (1,1,0)\otimes (0,0,1) &=(1,1,1)\cup 2(2,1,0)\cup 2(1,0,1) \cup 4(0,1,0)\,.
\intertext{}
H_3\ :\quad (1,0,0)\otimes (0,0,1) &=(1,0,1)\cup (0,\tau-1,\tau-1)\cup 5(\tau,0,0) \cup\\
                                \notag & \quad 3(0,0,\tau-1)\,.
\end{align*}

\subsection{Decomposition of products of $E_8$ orbits.}\

We say that an orbit is fundamental if its dominant weight in
the $\o$-basis has precisely one coordinate equal to 1 and all others
are zero. Thus $E_8$ has 8 fundamental orbits. Their sizes
range from 240 to over 17\,000.

All 36 different products of fundamental orbits of $E_8$ were
decomposed in \cite{GP} and are explicitly shown within the tables.
They were indispensable in solving the main problem of
\cite{GP}, namely the decomposition of products of fundamental
representations of $E_8$.

\section{Decomposition of symmetrized powers of orbits}

\subsection{Symmetrized powers of $G$-polynomials}\

The product of $m$ identical polynomials, say $P(\lambda;x)$, is the
subject of the action of the permutation group $S_m$ of $m$
elements. Thus it can be decomposed into a sum of components with a
specific permutation symmetry. It is well known from representation
theory that the permutation symmetry commutes with the action of the
Weyl group. Consequently, each permutation symmetry component can be
decomposed into a sum of polynomials.

Let ${\tiny\yng(1)}$ be short-hand notation for a polynomial
\eqref{polynomial}. The product of two and more copies of
${\tiny\yng(1)}$ decomposes into the symmetry components indicated
by their Young tableaux:
\begin{equation}
{\tiny\yng(1)}\otimes{\tiny\yng(1)} = {\tiny\yng(1,1)}
     +{\tiny\yng(2)}\,,
     \qquad
{\tiny\yng(1)}\otimes{\tiny\yng(1)}\otimes{\tiny\yng(1)} =
{\tiny\yng(1,1,1)}+2\,{\tiny\yng(2,1)}+ {\tiny\yng(3)}\,,\ \ \dots\
\end{equation}

In general, the square stands for a set of $G$-invariant items, each
square containing the same items. Those can be monomials of a
polynomial, or weights in the case of the weight system of a
representation of a semisimple Lie group/algebra, or points of a
$G$-orbit. The product of $m$ copies of the same square decomposes
into permutational symmetry components according to the
representations of the group $S_m$. The components are identified by
their Young tableau.  Each of the components is further decomposable
into the sum of parts that are labeled by the orbits of the Coxeter
group $G$.

In order to perform such a two-step decomposition, (i) the items of
the square need to be numbered consecutively in any convenient way.
The items belonging to a particular permutation symmetry component
are then determined according to the inequalities, shown in the next
subsection, and more generally implied by the corresponding Young
tableau. Then (ii) items belonging to a particular Young tableau,
which are labeled by points transformed by $G$, are sorted out into
the Coxeter group orbits. Practically it suffices to find the items
labeled by dominant points.

\subsection{Symmetrized powers of $G$-orbits}\

For simplicity of notation let us continue to label an orbit
$G(\lambda)$ by its dominant point $\lambda$. The product of the
same two $G$-orbits decomposes into its symmetric and antisymmetric
parts:
\begin{equation}
\lambda\otimes\lambda=(\lambda^2)_{symm}
                     \cup(\lambda^2)_{anti}
\end{equation}

Each of the two terms of the right side is further decomposable into
the sum of individual orbits. Let $\lambda_1,\lambda_2,\dots$ be the
points of the orbit $\lambda$ numbered in any order. Then the
content of the two parts is determined by the following inequalities
\begin{alignat}{2}
&(\lambda^2)_{symm}
     &&\ni \lambda_p+\lambda_q\,,\quad p\geq q\,,\\
&(\lambda^2)_{anti}
     &&\ni \lambda_p+\lambda_q\,,\quad p>q\,.
\end{alignat}
The product of 3 copies of $\lambda$ decomposes likewise
\begin{equation}
\lambda\otimes\lambda\otimes\lambda
      =(\lambda^3)_{symm}\cup(\lambda^3)_{anti}\cup2(\lambda^3)_{mixed}\,,
\end{equation}
where permutation symmetry components are formed from the $N$
points as follows:
\begin{alignat}{2}
&(\lambda^3)_{symm}&&\ni\lambda_p+\lambda_q+\lambda_s\,,\quad
     p\geq q\geq s\,,\\
&(\lambda^3)_{anti}&&\ni\lambda_p+\lambda_q+\lambda_s\,,\quad
     p>q>s\,,\\
&(\lambda^3)_{mixed}&&\ni\lambda_p+\lambda_q+\lambda_s\,,\quad
     p\geq q\ \text{and}\  p>s\,.
\end{alignat}

Similarly, any higher power decomposes into permutation symmetry
components where each is a sum of individual orbits.

\subsection{Two-dimensional examples}\
\begin{alignat*}{2}
A_2\ :\quad &(0,1)^2_{symm} &&=(1,0)\cup(0,2)\,,\\
            &(0,1)^2_{anti} &&=(1,0)\,.\\
            &(1,1)^2_{symm} &&=(2,2)\cup(1,1)\cup(3,0)\cup(0,3)\cup3(0,0)\,,\\
            &(1,1)^2_{anti} &&=(1,1)\cup(3,0)\cup(0,3)\cup3(0,0)\,.\\
            &(1,0)^3_{symm} &&=(1,1)\cup(3,0)\cup(0,0)\,,\\
            &(1,0)^3_{anti} &&=(0,0)\,,\\
            &(1,0)^3_{mixed}&&=(1,1)\cup2(0,0)\,.
\intertext{}
C_2\ :\quad &(0,1)^2_{symm} &&=(2,0)\cup(0,2)\cup2(0,0)\,,\\
            &(0,1)^2_{anti} &&=(2,0)\cup2(0,0)\,.\\
            &(1,0)^3_{symm} &&=(1,1)\cup(3,0)\cup2(1,0)\,,\\
            &(1,0)^3_{anti} &&=(1,0)\,,\\
            &(1,0)^3_{mixed}&&=(1,1)\cup3(1,0)\,.\\
\intertext{}
G_2\ :\quad &(0,1)^2_{symm} &&=(1,0)\cup(0,2)\cup(0,1)\cup3(0,0)\,,\\
            &(0,1)^2_{anti} &&=(1,0)\cup(0,1)\cup3(0,0)\,.\\
            &(1,0)^3_{symm} &&=(1,3)\cup(3,0)\cup(2,0)\cup3(1,0)\cup2(0,3)\cup2(0,0)\,,\\
            &(1,0)^3_{anti} &&=(2,0)\cup2(1,0)\cup2(0,0)\,,\\
            &(1,0)^3_{mixed}&&=(1,3)\cup2(2,0)\cup5(1,0)\cup2(0,3)\cup4(0,0)\,.\\
\intertext{}
H_2\ :\quad &(0,1)^2_{symm} &&=(\tau,0)\cup(0,2)\cup(0,\tau-1)\,,\\
            &(0,1)^2_{anti} &&=(\tau,0)\cup(0,\tau-1)\,.\\
            &(1,0)^3_{symm} &&=(2-\tau,1)\cup(1,\tau)\cup(3,0)\cup(\tau,0)\cup(0,\tau-1)\,,\\
            &(1,0)^3_{anti} &&=(\tau,0)\cup(0,\tau-1)\,,\\
            &(1,0)^3_{mixed}&&=(2-\tau,1)\cup(1,\tau)\cup2(\tau,0)\cup2(0,\tau-1)\,.\\
\end{alignat*}

\subsection{Three-dimensional examples}\
\begin{alignat*}{2}
A_3\ :\quad &(1,0,0)^3_{symm} &&=(1,1,0)\cup(3,0,0)\cup(0,0,1)\,,\\
            &(1,0,0)^3_{anti} &&=(0,0,1)\,,\\
            &(1,0,0)^3_{mixed}&&=(1,1,0)\cup2(0,0,1)\,.\\
\intertext{}
B_3\ :\quad &(1,0,0)^3_{symm} &&=(1,1,0)\cup(3,0,0)\cup3(1,0,0)\cup(0,0,2)\,,\\
            &(1,0,0)^3_{anti} &&=2(1,0,0)\cup(0,0,2)\,,\\
            &(1,0,0)^3_{mixed}&&=(1,1,0)\cup5(1,0,0)\cup2(0,0,2)\,.\\
\intertext{}
C_3\ :\quad &(1,0,0)^2_{symm} &&=(2,0,0)\cup(0,1,0)\cup3(0,0,0)\,,\\
            &(1,0,0)^2_{anti} &&=(0,1,0)\cup3(0,0,0)\,.\\
\intertext{}
H_3\ :\quad &(1,0,0)^2_{symm} &&=(2,0,0)\cup(0,1,0)\cup(0,\tau-1,0)\cup6(0,0,0)\,,\\
            &(1,0,0)^2_{anti} &&=(0,1,0)\cup(0,\tau-1,0)\cup6(0,0,0)\,.\\
\end{alignat*}


\section{Congruence classes, indices, and anomaly numbers of polytopes}

Here we introduce numerical characterizations of $W$-orbits, analogs
of similar quantities known for irreducible representations of
semisimple Lie groups, which proved particularly useful in their application.

\subsection{Congruence classes}\

Inclusion among the lattices \eqref{lattices} is an important
property of the Weyl group $W$. The weight lattice $P$ can be
understood as a union of several components, each isomorphic to the
root lattice $Q$. The components are shifted relative to each other
by some elements of $P$. An individual component consists of points
belonging to one congruence class of $P$. The index of $Q$ in $P$,
denoted $|Z|$, is the number of distinct congruence classes in $P$.
The value of $|Z|$ reflects other properties of $G$. For example, it is
the order of the center of $G$, it is a common denominator of
coordinates of all points of $P$ when given in the basis of simple
roots, etc. One has $|Z|>1$ for all $G$ but for the exceptional
simple Lie groups of types $E_8$, $F_4$, and $G_2$.

The congruence number $c$ is a number attached to points of $P$. The
value of $c$ is common to all points of the same congruence class.
It can be defined in a number of equivalent ways. Our definition
coincides with that of \cite{LP}. All points of any $W$-orbit belong
to the same congruence class. Furthermore, orbits obtained from the
decomposition of a product belong to the same congruence class, and
their congruence number is the sum of the congruence numbers of the
orbits of the multiplication. That is also true for the
decomposition of symmetrized powers of orbits.

Let $x=(x_1,x_2,\dots,x_n)\in P$ be a point to consider in the $\o$-basis. Its congruence number $c(x)$ is given by the following formulas :
\begin{equation}
\begin{aligned}
A_n\quad :\quad c(x)&=\sum_{k=1}^n kx_k\mod(n+1)\\
B_n\quad :\quad c(x)&=x_n\mod2\\
C_n\quad :\quad c(x)&=\sum_{k=1}^{[\tfrac{n+1}2]} x_{2k-1}\mod2\\
D_n\quad :\quad c(x)&=(c_1(x) \mod2\,,\ c_2(x)\mod4),\\
                 c_1(x)&=x_{n-1}+x_n\\
                 c_2(x)&=\begin{cases}
                     2x_1+2x_3+\dots+2x_{n-2}+(n-2)x_{n-1}
                                       +nx_n,\quad n\ \text{odd}\\
                     2x_1+2x_3+\dots+2x_{n-3}+(n-2)x_{n-1}+nx_n,\quad n\
                          \text{even}
\end{cases}\\
E_6\quad :\quad c(x)&=x_1 - x_2 + x_4 - x_5 \mod3\\
E_7\quad :\quad c(x)&=x_4 + x_6 + x_7 \mod2
\end{aligned}
\end{equation}

For $E_8, F_4$ and $G_2$ there is only one congruence class, namely
$c(x)=0$ for all $x \in P$. Note also that the roots of any group
belong to the congruence class $c(x)=0$. Hence also the points of
the root lattice of any group belong to the congruence class
$c(x)=0$.

The points of any single $G$-orbit belong to the same congruence
class because the difference between any two points of the same
orbit is an integer linear combination of simple roots, as can be
derived from \eqref{michelle}.

For the non-crystallographic groups, the congruence classes can be similarly
defined, involving their appropriate irrationality.
It is important to recall that, in these cases, $P$ is a dense
lattice. The coordinates of $x\in P$, relative to the $\o$-basis, are
the numbers $a+\tau b$, with $a,b\in\Z$.
\begin{equation}
H_2\quad : \quad c(x)= \tau x_1 + 2x_2 \mod5, \quad \text{where $\tau=3$}\\
\end{equation}

\subsection{The second and higher indices}\

The second and higher indices were defined~\cite{PSW} for weight
systems of irreducible finite dimensional representations of compact
semisimple Lie groups. Extensive tables of indices of degree 0, 2
and 4 are found in~\cite{McP}. The fact that a weight system is a
union of several $W$-orbits suggests that the indices could be
introduced for individual orbits. Moreover, we introduce them also
for non-crystallographic Coxeter groups with the same formulas.

For any finite Coxeter group $G$, we define an index
$I_\lambda^{(2k)}$ of degree $2k$ of a $G$-orbit $G(\lambda)$ by
\begin{equation}\label{index}
I_\lambda^{(2k)}=\sum_{\mu\in G(\lambda)}\l\mu,\mu\r^k
  =\l\lambda,\lambda\r^kI_\lambda^{(0)}\,,\qquad
        k=0,1,2,\dots\,,
\end{equation}
because points of $G(\lambda)$ are equidistant from the origin.
Clearly $I_\lambda^{(0)}=|G(\lambda)|$ is the number of points of
the orbit $G(\lambda)$ given by \eqref{orbitsize}.

Higher indices of products of two orbits, $G(\lambda_1)\otimes
G(\lambda_2)$, are also useful in calculating the decompositions.
Let $r$ be the rank of $G$.

\begin{equation}\label{index2}
I^{(2k)}(G(\lambda_1)\otimes G(\lambda_2))=I_{\lambda_1\otimes \lambda_2}^{(2k)}=I_{\lambda_1 + \lambda_2}^{(2k)} + \dots + I_{\lambda_1 + \overline{\lambda_2}}^{(2k)}\\
\end{equation}

\begin{align}
I_{\lambda_1\otimes \lambda_2}^{(0)}
     &=I_{\lambda_1}^{(0)}\,I_{\lambda_2}^{(0)}\\
I_{\lambda_1\otimes \lambda_2}^{(2)}
     &=I_{\lambda_1}^{(2)}\,I_{\lambda_2}^{(0)}
     +I_{\lambda_1}^{(0)}\,I_{\lambda_2}^{(2)}\label{I2} \\
     &=I_{\lambda_1}^{(0)}\,I_{\lambda_2}^{(0)}\left(
     \l\lambda_1,\lambda_1\r+\l\lambda_2,\lambda_2\r\right)\\
I_{\lambda_1\otimes \lambda_2}^{(4)}
     &=I_{\lambda_1}^{(4)}\,I_{\lambda_2}^{(0)}
     +\frac{2(r+2)}{r}I_{\lambda_1}^{(2)}\,I_{\lambda_2}^{(2)}
     +I_{\lambda_1}^{(0)}\,I_{\lambda_2}^{(4)}
\end{align}

Table 4 presents examples of indices of degree 0, 2, 4, 6 and 8 for individual orbits of $A_2$, $C_2$, $G_2$ and $H_2$. 

\begin{table}
\begin{center}
\begin{tabular}{|c|c|c|c|c|c|}
\hline
$A_2$&$I^{(0)}$ & $I^{(2)}$ & $3I^{(4)}$ & $9I^{(6)}$ & $27I^{(8)}$ \\
\hline
$(1,0)$&3 &2 &4 &8 &16\\
\hline
$(2,0)$&3 &8 &64 &512 &4096\\
\hline
$(1,1)$&6 &12 &72 &432 &2592\\
\hline
$(2,1)$&6 &28 & 392 &5488 &76832\\
\hline
\end{tabular}
\medskip

\begin{tabular}{|c|c|c|c|c|c|}
\hline
$C_2$&$I^{(0)}$ & $I^{(2)}$ & $2I^{(4)}$ & $4I^{(6)}$ & $8I^{(8)}$ \\
\hline
$(1,0)$& 4 & 2& 2 &2 &2\\
\hline
$(0,1)$& 4 & 4 & 8 &16 &32 \\
\hline
$(2,0)$& 4 & 8 & 32 & 128 &512 \\
\hline
$(0,2)$& 4 & 16 & 128 &1024 &8192 \\
\hline
$(1,1)$& 8 & 20 & 100 &500 &2500 \\
\hline
$(2,1)$& 8 & 40 & 400 &4000 &40000 \\
\hline
\end{tabular}
\medskip

\begin{tabular}{|c|c|c|c|c|c|}
\hline
$G_2$&$I^{(0)}$ & $I^{(2)}$ & $3I^{(4)}$ & $9I^{(6)}$ & $27I^{(8)}$\\
\hline
$(0,1)$&6 & 4 & 8 &16 &32  \\
\hline
$(1,0)$&6 &12 &72 &432 &2592 \\
\hline
$(0,2)$&6 & 16 & 128 &1024 &8192 \\
\hline
$(0,3)$&6 & 36 & 648 &11664 &209952 \\
\hline
$(2,0)$&6 & 48 & 1152 &27648 &663552 \\
\hline
$(1,1)$&12 & 56 & 784 &10976 &153664 \\
\hline
\end{tabular}
\medskip

\begin{tabular}{|c|c|c|c|}
\hline
$H_2$&$I^{(0)}$ & $(3-\tau)I^{(2)}$ & $(3-\tau)^2I^{(4)}$\\
\hline
$(1,0)$& 5 & 10 & 20\\
\hline
$(2,0)$& 5 & 40 & 320\\
\hline
$(1,1)$& 10 & $20(\tau+2)$ & $40(\tau+2)^2$\\
\hline
$(2,1)$& 10 & $10(4\tau+10)$ & $10(4\tau+10)^2$\\
\hline
\end{tabular}
\medskip
\caption{Examples of the indices $I^{(2k)}, k = 0,\dots,4$.}
\label{tab4}
\end{center}
\end{table}

\subsection{Anomaly numbers}\

Triangle anomaly numbers were introduced in physics \cite{S,GJ,O8}
as quantities assigned to irreducible representations of a few
compact semisimple Lie groups and calculated from the weight systems
of their representations. Constraints on possible models in particle
physics were imposed in terms of admissible values of the anomaly
numbers of representations involved in a particular model.
Generalization of the concept to all compact semisimple Lie groups
and to higher than third degree anomaly number originates in
\cite{PS8}. Our goal here is to show that the anomaly numbers can be
used also for constituents of the weight systems of irreducible
representations, namely for $W$-orbits and more generally, for the
orbits $G(\lambda)$ of any finite Coxeter group.

The anomaly number $I_\lambda^{(2k-1)}$ of degree $2k-1$ of the
orbit $G(\lambda)$ of the Coxeter group $G$ is defined as follows,
\begin{equation}\label{anomaly}
I_\lambda^{(2k-1)}=\sum_{\mu\in
G(\lambda)}\l\mu,u\r^{2k-1}\,,\qquad k=1,2,\dots\,,
\end{equation}
where $u$ is a special vector passing through the origin. In
particular, $I^{(1)}=0$ in all cases. The anomaly number of physics
literature is $I^{(3)}$, therefore it is the only one we consider.

Frequently used property of $I^{(3)}$ is the decomposition of the
product of two orbits, which is the analog of \eqref{I2}:
\begin{equation}
I_{\lambda_1\otimes \lambda_2}^{(3)}
     =I_{\lambda_1}^{(3)}\,I_{\lambda_2}^{(0)}
     +I_{\lambda_1}^{(0)}\,I_{\lambda_2}^{(3)}
     =I_{\lambda_1 + \lambda_2}^{(3)} + \dots + I_{\lambda_1 + \overline{\lambda_2}}^{(3)}
\end{equation}

In general terms, the direction of $u$ can be characterized as
follows. Suppose $W$ in \eqref{anomaly} is the Weyl group of a
compact simple Lie group $G$, and that $G$ has a maximal reductive
subgroup of type $U(1)\times G'$. Then the direction of $u$ is given
by the direction corresponding to $U(1)$ in the Euclidean space
spanned by the roots of $G$.

The first question to answer is when such a maximal subgroup is
present. For a complete list of the cases see below \cite{BdeS}:
\begin{equation}
\begin{alignedat}{2}
A_n &\supset A_{n-1}\times U(1) \qquad &&n\geq2\\
A_n &\supset A_k\times A_{n-k-1}\times U(1) \qquad
               &&n\geq3\,,\quad 1\leq k\leq[\tfrac{n-1}2]\\
B_n &\supset B_{n-1}\times U(1) \qquad &&n\geq3\\
C_n &\supset A_{n-1}\times U(1) \qquad &&n\geq2\\
D_n &\supset A_{n-1}\times U(1) \qquad &&n\geq4\\
D_n &\supset D_{n-1}\times U(1) \qquad &&n\geq5\\
E_6 &\supset D_5\times U(1) \\
E_7 &\supset E_6\times U(1)
\end{alignedat}
\end{equation}
As long as each orbit of a given group contains with every weight
also its negative, the anomaly numbers are equal to zero. Therefore
the interesting cases that remain are found in $A_n$, $D_{2k+1}$,
$E_6$, and $E_7$. In physics, however, only the anomaly numbers of
$A_n\supset A_{n-1}\times U(1)$ are used so far.

Anomaly numbers of $H_2$, $H_3$, and $H_4$ are also defined by
\eqref{anomaly}. In those cases, however, the direction of $u$ has
to be determined differently since there is no $U(1)$ subgroup.
Instead, one can require that $u$ be orthogonal to selected simple
roots: $\alpha_1$ for $H_2$,   $\alpha_1$ and $\alpha_2$ for $H_3$,
and $\alpha_1$, $\alpha_2$, and $\alpha_3$ for $H_4$. Anomaly
numbers for $H_2$ are zero for all orbits. They will be considered
elsewhere \cite{HLP}, along with the anomaly numbers of other
non-crystallographic groups.

 \section{Concluding remarks}
\begin{enumerate}

\item
Useful and interesting objects may turn out to be $G$-orbits with
each point decorated by a sign \cite{P} according to the following
rule. The dominant point, say $\lambda$, and all points obtained
from it by an even number of reflections generating $G$, carry a
positive sign, while all points of the orbit obtained from
$\lambda$ by an odd number of reflections carry a negative sign. Let
us call an $S$-orbit a decorated orbit of $\lambda$ of $G$, while the
orbits without the sign decoration, i.e. all positive signs, are
called $C$-orbits of $\lambda$ of $G$.  In order to avoid ambiguities,
it should be stipulated that $\lambda$ of an $S$-orbit must have all
coordinates positive in $\o$-basis.

Multiplication of such orbits follows simple rules:
\begin{alignat}{2}
\text{$C$-orbit} &\times\text{$C$-orbit}\ &&\longrightarrow\
\text{$C$-orbits},
\label{CC}\\
\text{$C$-orbit} &\times\text{$S$-orbit}\ &&\longrightarrow\ \text{$S$-orbits},          \label{CS}\\
\text{$S$-orbit} &\times\text{$S$-orbit}\ &&\longrightarrow\
\text{$C$-orbits}. \label{SS}
\end{alignat}
In \eqref{CC}, all coefficients in the decomposition of the product
are positive integers, while in \eqref{CS} and \eqref{SS}, all such
coefficients are integers, but not all may be positive.

The decomposition of many products of $C$-orbits with lowest nontrivial
$S$-orbit can be directly inferred from the tables \cite{BMP}, using the
Weyl character formula.
\smallskip

\item
In the examples, we often required that a $G$-orbit consist of
points of the weight lattice $P$. Very few properties of the orbits
would have been lost, had we instead allowed $\lambda\in\R^n$. The
congruence classes would not then be applicable.

Consider the following products of $A_2$ orbits as examples:
\begin{alignat}{3}
(a,0)&\otimes(\varepsilon,0)
    &&=(a+\varepsilon,0)\cup(a-\varepsilon,\varepsilon)\,,
    \qquad &&0<\varepsilon\ll 1\,,\ a\geq1\,,\label{close}\\
(a,0)&\otimes(0,a+\varepsilon)
    &&=(a,a+\varepsilon)\cup(0,\varepsilon)\,,
    \quad &&0<\varepsilon\ll 1\,,\ a\gg1\label{thiny}\,.
\end{alignat}
The radii of the two orbits in the decomposition \eqref{close} can
be drawn arbitrarily close by a suitable choice of $\varepsilon$, and
in \eqref{thiny} they can be pushed as far apart as desired by the choice of $a$. The
second orbit in \eqref{thiny} has a radius equal to
$\varepsilon\sqrt{\tfrac23}$.
\smallskip

\item
For a geometric interpretation of orbits as polytopes, refer to the
paragraph following equation \eqref{decomposition}. The `interaction' (i.e.
product) between two concentric orbit-layers results in the layered
structure of orbits. They are subject to the equality of indices of
various degrees, congruence numbers, relations between anomaly
numbers. Speculative interpretation can go further: Consider
$I_\lambda^{(2)}$ as the `energy' of the orbit and
$I_{\lambda\otimes\lambda'}^{(2)}$ as the `energy' of the
interacting pair, etc.
\smallskip

\item
Although we did not pursue it here, orbit
multiplication can be viewed as an `interaction' between two orbits similarly as
used in particle physics to view interacting multiplets of
particles. A multiplet is described by the weight system of an
irreducible representation of the corresponding Lie group/algebra.
Here, the role of the multiplet would be given to the set of points
of an orbit. In both cases, such interactions would be governed by the
strict equality of indices of various degrees, congruence numbers,
relations between anomaly numbers.
But there is a price to pay for such a reinterpretation of
multiplets: the overall invariance of the theory with respect to the
Lie group would be reduced to the invariance with respect to the
Coxeter group, or to its (discrete) image `lifted' into the Lie
group \cite{MP.DT}.
\smallskip

\item
It would be useful to ask additional questions about the properties
of indices and anomaly numbers of various degrees. Such questions
can be answered by adaptation of the methods used for the weight
system of representations \cite{PSW, PS8}.
\smallskip

\item
In place of finite Coxeter groups, we could have chosen to consider
other finite groups for similar considerations \cite{Macd}.
The immediate motivations for our choice were recent applications in
harmonic analysis, where $W$-orbits are playing a fundamental
role. Equally interesting would be to consider orbits of infinite
Coxeter groups. (A Coxeter group with connected diagram is of
infinite order if its diagram is different from those listed in
Section~2.) The orbits of representations of Kac-Moody algebras would be relatively easily amenable to such a study.
\smallskip

\item
Similarly, we could consider orbits of two or more seed points. A simple
example is the root system of the group $G_2$. Choosing as the two seed points one short root and one long root, say
$\alpha_2$ and $\alpha_1+3\alpha_2$, the orbit of the pair is a star-like polygon formed by the root system of $G_2$.
\smallskip

\item
An interesting problem appears to be to pursue a similar study of
orbits of the even subgroups of Coxeter groups, particularly because
these subgroups are not Coxeter groups in general.
\end{enumerate}

\bigskip
\centerline{\bf Acknowledgements}
\medskip

Work supported in part by the Natural Sciences and Engineering
Research Council of Canada, the MIND Research Institute.
We are grateful for the hospitality extended
to us by the Doppler Institute, Czech Technical
University (J.P. and M.L.), and by the Centre de Recherches
Math\'ematiques, Universit\'e de Montr\'eal (L.H.).


\end{document}